\documentclass[A4]{article}
\usepackage{amscd, amssymb, amsmath, amsthm}
\usepackage{latexsym}
\usepackage{upref}
\usepackage{epsfig}
\usepackage{mathrsfs}
\usepackage{float, floatflt}
\usepackage{indentfirst}
\usepackage{hyperref}
\hypersetup{
    colorlinks=true,
    linkcolor=blue,
    citecolor=blue,
    filecolor=blue,
    urlcolor=blue,
}
\usepackage{graphics,graphicx,color}
\usepackage{cite}
\usepackage{enumitem}
\usepackage{fancyhdr} 
\usepackage[utf8]{inputenc}         
\usepackage[T1]{fontenc}
\usepackage{mathptmx} 

\AtBeginDocument{
  \DeclareSymbolFont{AMSb}{U}{msb}{m}{n}
  \DeclareSymbolFontAlphabet{\mathbb}{AMSb}}  

\usepackage{rotating} 
\usepackage{booktabs} 
\usepackage{pb-diagram,pb-xy}       
\usepackage{tikz}                   
\usetikzlibrary{matrix,arrows}
\usepackage{array, caption} 

\theoremstyle{plain}
\newtheorem{theorem}{Theorem}[section]

\newtheorem{proposition}[theorem]{Proposition}

\newtheorem{corollary}[theorem]{Corollary}
\theoremstyle{definition}
\newtheorem{definition}[theorem]{Definition}
\newtheorem{example}[theorem]{Example}
\newtheorem*{question}{Question}
\theoremstyle{remark}

\numberwithin{equation}{section}

\makeatletter
\def\th@plain{%
  \thm@notefont{}
  \itshape 
}
\def\th@definition{%
  \thm@notefont{}
  \normalfont 
} \makeatother

\setlist{font=\normalfont}

\DeclareMathAlphabet{\cols}{OMS}{cmsy}{m}{n} %
\DeclareMathAlphabet{\mathcal}{OMS}{cmsy}{m}{n}

\newcommand{\If}{\Rightarrow}
\newcommand{\Fi}{\Leftarrow}
\newcommand{\set}[1]{\{#1\}}
\newcommand{\cset}[2]{\set{{#1}\colon{#2}}}
\newcommand{\lcup}[2]{\displaystyle\bigcup_{#1}^{#2}}
\newcommand{\abs}[1]{|#1|}
\newcommand{\gap}{\hskip0.3em}
\newcommand{\Fix}[1]{\mathrm{Fix}\,{#1}}
\newcommand{\sym}[1]{\mathrm{Sym}\,{(#1)}}
\newcommand{\stab}[1]{\mathrm{Stab}\,{(#1)}}
\newcommand{\gyr}[2]{{\mathrm{gyr}[{#1}]}{#2}}
\newcommand{\aut}[1]{\mathrm{Aut}\,{(#1)}}
\newcommand{\GYR}[1]{\mathrm{GYR}\,{(#1)}}
\newcommand{\gen}[1]{\langle#1\rangle}
\newcommand{\igyr}[2]{{\mathrm{gyr^{-1}}[{#1}]}{#2}}
\newcommand{\qt}[1]{``#1''}
\newcommand{\res}[2]{{#1}\hskip-1pt\mid_{#2}}
\newcommand{\D}{\mathbb{D}}
\newcommand{\vphi}{\varphi}
\newcommand{\norm}[1]{\|#1\|}
\newcommand{\B}{\mathbb{B}}
\newcommand{\R}{\mathbb{R}}
\newcommand{\Bp}[1]{\left(#1\right)}
\newcommand{\N}{\mathbb{N}}

\newcommand{\C}{\mathbb{C}}
\renewcommand{\vec}[1]{\mathbf{#1}}

\begin{document}
\title{\textbf{Geometry of gyrogroups via Klein's approach}}
\author{Teerapong Suksumran\,\href{https://orcid.org/0000-0002-1239-5586}{\includegraphics[scale=1]{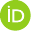}}\\
Research Center in Mathematics and Applied Mathematics\\
Department of Mathematics\\
Faculty of Science, Chiang Mai University\\
Chiang Mai 50200, Thailand\\
\small{E-mail address: {\tt teerapong.suksumran@cmu.ac.th}}}
\date{}
\maketitle

\begin{abstract}
Using Klein's approach, geometry can be studied in terms of a space of points and a group of transformations of that space. This allows us to apply algebraic tools in studying geometry of mathematical structures. In this article, we follow Klein's approach to study the geometry $(G, \mathcal{T})$, where $G$ is an abstract gyrogroup and $\mathcal{T}$ is an appropriate group of transformations containing all gyroautomorphisms of $G$. We focus on $n$-transitivity of gyrogroups and also give a few characterizations of coset spaces to be minimally invariant sets. We then prove that the collection of open balls of equal radius is a minimally invariant set of the geometry $(G, \Gamma_m)$ for any normed gyrogroup $G$, where $\Gamma_m$ is a suitable group of isometries of $G$.
\end{abstract}
\thispagestyle{empty}

\section{Introduction}
Geometry via Klein's approach, the so-called {\it Erlangen Program}, can be studied by specification of a set of points and a group of transformations. Using this approach, \mbox{geometry} is the study of objects and properties that remain invariant under allowable transformations. One of the most advantages is that three well-known geometries (Euclidean, hyperbolic, and elliptic geometries) emerge as special cases of general geometries. For instance, Euclidean geometry may be studied via the space $\R^n$ of $n$-tuples with real components, together with the group of transformations given by
$$
\cols{E}(\R^n) = \cset{L_\vec{u}\circ \tau}{\vec{u}\in\R^n, \tau\in\mathrm{O}(\R^n)},
$$
where $L_\vec{u}$ is the left translation $L_{\vec{u}}\colon \vec{v}\mapsto \vec{u}+\vec{v}$ and $\mathrm{O}(\R^n)$ is the group of orthogonal transformations on $\R^n$. Moreover, hyperbolic geometry may be studied via the set $\D$ of complex numbers with modulus less than $1$, together with the group of M\"{o}bius transformations on $\D$ given by
$$
\cols{M}(\D) = \cset{\tau_a\circ \rho}{a\in\D, \rho\in\mathrm{O}(\D)},
$$
where $\tau_a$ is the M\"{o}bius translation $\tau_a\colon z \mapsto \dfrac{a+z}{1+\bar{a}z}$ and $\mathrm{O}(\D)$ is the group of disk rotations $\rho\colon z\mapsto \omega z$ with $\abs{\omega} = 1$. Klein's geometry, also known as transformation geometry, is the philosophy behind the notion of a symmetry group. A nice application of this approach can be found in a study of physics from symmetry; see, for instance, \cite{MR3751325}.

The notion of a gyrogroup was introduced by A. Ungar in the study of parametrization of the Lorentz transformation group \cite{AU1988TRP}. Since then gyrogroups have been widely studied in several fields. In particular, analytic geometry of concrete gyrogroups such as Einstein and M\"{o}bius gyrogroups has been intensively examined by using a gyrovector space approach in a series of books; see, for instance, \cite{MR3308966, MR2723718, MR2919420, AU2008AHG}. From the algebraic viewpoint, gyrogroups may be regarded as algebraic structures that generalize groups. In fact, they share common properties with groups. Recently, many algebraic properties of abstract gyrogroups have been released. This combined with the two aforementioned examples inspires us to study geometry of arbitrary gyrogroups using Klein's approach in the present article. The results here give us insight into both algebraic and geometric properties of gyrogroups. 

\section{Preliminaries}
\subsection{Basic knowledge of Klein's geometry}
In this section, we review and prove some results in the basic theory of geometries using Klein's approach. This approach enables us to study geometry in terms of a space of points and a group of transformations of that space. See, for instance, \cite{hitchman2018geometry}. This will prove useful in studying geometry of abstract gyrogroups in Sections \ref{sec: Geometry of gyrogroups} and \ref{sec: normed gyrogroups}.

\begin{definition}
A geometry $(S, \cols{T})$ consists of a nonempty set $S$, together with a group $\cols{T}$ of bijections from $S$ to itself under function composition $\circ$.
\end{definition}

Let $(S, \cols{T})$ be a geometry. Any element of $S$ is called a {\it point}; any subset of $S$ is called a {\it figure}; any bijection in $\cols{T}$ is called a {\it transformation} on $S$. A set $\cols{F}$ of figures in $S$ is {\it invariant} if $T(F)$ belongs to $\cols{F}$ for all $F\in\cols{F},T\in\cols{T}$. An invariant set is {\it minimal} if it contains no nonempty proper invariant subset. A mathematical property related to the geometry $(S, \cols{T})$ is {\it $\cols{T}$-invariant} (or simply {\it invariant} if there is no confusion) if it is preserved under any transformation in $\cols{T}$. Denote by $\cols{P}(S)$ the power set of $S$. A function $f$ defined on $\cols{P}(S)$ is {\it invariant} if $f(A) = f(T(A))$ for all $A\in\cols{P}(S), T\in\cols{T}$.

The familiar notion of congruence in Euclidean geometry can be stated in terms of transformations as follows.

\begin{definition}
Let $(S, \cols{T})$ be a geometry. We say that figures $A$ and $B$ in $S$ are {\it congruent}, written $A\cong B$, if there exists a transformation $T$ in $\cols{T}$ such that $T(A) = B$.
\end{definition}

In any geometry $(S, \cols{T})$, as $\cols{T}$ is a group under composition, we immediately obtain that $\cong$ is an equivalence relation on the power set of $S$. The next theorem shows a strong connection between the minimally invariant sets and the equivalence classes determined by the congruence relation.

\begin{theorem}\label{thm: characterization of minimally invariant set}
Let $(S, \cols{T})$ be a geometry and let $\cols{O}(S, \cols{T})$ be the collection of equivalence classes determined by the congruence relation associated to $(S, \cols{T})$. Then $\cols{C}\in \cols{O}(S, \cols{T})$ if and only if $\cols{C}$ is a nonempty minimally invariant set.
\end{theorem}
\begin{proof}
($\If$) Let $\cols{C}\in \cols{O}(S, \cols{T})$. By definition, $\cols{C}\ne\emptyset$. Let $A\in\cols{C}$ and let $T\in\cols{T}$. Since $A\cong T(A)$ and $A\in\cols{C}$, it follows that $T(A)\in\cols{C}$. This shows that $\cols{C}$ is an invariant set. Suppose that $\cols{D}$ is a nonempty proper subset of $\cols{C}$. Then there is a figure $B$ such that $B\in\cols{C}$ but $B\not\in\cols{D}$. Since $\cols{D}\ne\emptyset$, there is a figure $A\in\cols{D}$. Since both $A$ and $B$ lie in $\cols{C}$, we have $A\cong B$ and so $B = T(A)$ for some $T\in\cols{T}$. Hence, $A\in\cols{D}$ but $T(A)\not\in\cols{D}$, which shows that $\cols{D}$ is not invariant. Thus, $\cols{C}$ is minimally invariant.

($\Fi$) Suppose that $\cols{C}$ is a nonempty minimally invariant set. Since $\cols{C}\ne\emptyset$, there is a figure $A\in\cols{C}$. Let $[A]$ be the equivalence class of $A$ determined by the congruence relation. By definition, $[A] = \cset{B\subseteq S}{B\cong A} = \cset{T(A)}{T\in\cols{T}}$. Furthermore, $[A]\in\cols{O}(S, \cols{T})$. We claim that $\cols{C} = [A]$. Since $\cols{C}$ is invariant and $A\in \cols{C}$, it follows that $[A]\subseteq \cols{C}$. Clearly, $[A]\ne\emptyset$. Next, we show that $[A]$ is invariant. Let $B\in [A]$
and let $T\in\cols{T}$. Then $B = S(A)$ for some $S\in\cols{T}$. Since $\cols{T}$ is a group, $T\circ S\in\cols{T}$. Thus, $T(B) = T(S(A)) = (T\circ S)(A)\in [A]$. This shows that $[A]$ is invariant. By definition, $[A]$ cannot be proper in $\cols{C}$. Hence, $[A] = \cols{C}$ and the proof completes.
\end{proof}

As an application of Theorem \ref{thm: characterization of minimally invariant set}, we obtain the following characterization of minimally invariant sets.

\begin{theorem}\label{thm: characterization of minimally invariant set II}
Suppose that $\cols{C}$ is a nonempty invariant set of a geometry $(S, \cols{T})$. Then $\cols{C}$ is minimally invariant if and only if $A\cong B$ for all $A, B\in\cols{C}$.
\end{theorem}

In light of Theorem \ref{thm: characterization of minimally invariant set}, we obtain a characterization of invariant sets related to the congruence relation as follows.

\begin{theorem}\label{thm: characterization of invariant set}
Let $(S, \cols{T})$ be a geometry and let $\cols{F}$ be a nonempty set of figures in $S$. Then $\cols{F}$ is an invariant set if and only if $\cols{F}$ is a union of equivalence classes determined by the congruence relation.
\end{theorem}
\begin{proof}
Suppose that $\cols{F}$ is an invariant set. Then $[A] = \cset{T(A)}{T\in\cols{T}}\subseteq \cols{F}$ for all $A\in\cols{F}$. It follows that $\cols{F} = \lcup{A\in\cols{F}}{}[A]$. Suppose conversely that $\cols{F} = \lcup{i\in I}{}[A_i]$, where $A_i$ is a figure in $S$ for all $i\in I$. Let $A\in\cols{F}$ and let $T\in\cols{T}$. Then $A\in [A_j]$ for some $j\in I$. Since $[A_j]$ is invariant, $T(A)\in[A_j]$ and so $T(A)\in \cols{F}$. This proves that $\cols{F}$ is invariant.
\end{proof}

The following theorem gives a characterization of invariant functions related to the congruence relation.

\begin{theorem}\label{thm: characterization of invariant function}
Let $(S, \cols{T})$ be a geometry, let $f$ be a function defined on $\cols{P}(S)$, and let $\cols{O}(S, \cols{T})$ be the collection of equivalence classes determined by the congruence relation associated to $(S, \cols{T})$. Then $f$ is invariant if and only if $f$ is constant on each class $\cols{C}$ in $\cols{O}(S, \cols{T})$.  
\end{theorem}
\begin{proof}
Suppose that $f$ is invariant. Let $\cols{C}\in\cols{O}(S, \cols{T})$ and let $X, Y\in \cols{C}$. Then $X\cong Y$ and so $Y = T(X)$ for some $T\in\cols{T}$. By assumption, $f(Y) = f(T(X)) = f(X)$. Thus, $f$ is constant on $\cols{C}$. Suppose conversely that $f$ is constant on $\cols{C}$ for all $\cols{C}\in \cols{O}(S, \cols{T})$. Let $A\in\cols{P}(S)$ and let $T\in\cols{T}$. By definition of congruence, $A\cong T(A)$ and so $A$ and $T(A)$ lie in the same class in $\cols{O}(S, \cols{T})$. By assumption, $f(A) = f(T(A))$. This shows that $f$ is invariant.
\end{proof}

\begin{corollary}
Let $(S, \cols{T})$ be a geometry and let $f$ be its invariant function. For all subsets $A, B$ of $S$, if $f(A)\ne f(B)$, then $A$ and $B$ are not congruent in $(S, \cols{T})$.
\end{corollary}

By definition, every transformation of a geometry is bijective. Hence, cardinality is an invariant function of any geometry.

\begin{theorem}
The cardinality is an invariant function of every geometry.
\end{theorem} 

\begin{corollary}
Let $(S, \cols{T})$ be a geometry. If $A$ and $B$ are subsets of $S$ such that $\abs{A}\ne\abs{B}$, then $A$ and $B$ are not congruent in $(S, \cols{T})$.
\end{corollary}

Usually, the study of geometries consists of
\begin{itemize}
\item the description of the congruence classes,
\item the identification of invariant properties, and
\item theorems about congruence classes and invariant properties.
\end{itemize}
According to Theorems \ref{thm: characterization of minimally invariant set}, \ref{thm: characterization of invariant set}, and \ref{thm: characterization of invariant function}, a complete description of minimally invariant sets, invariant sets, and invariant functions is obtained once the collection of equivalence classes determined by the congruence relation is known. Of course, it is difficult, in general, to determine all of equivalence classes associated to the congruence relation of a given geometry. 

Next, we relate the usual notion of homogeneity to the notion of invariance. \mbox{Recall} that a geometry $(S, \cols{T})$ is {\it homogeneous} if for each pair of points $x, y\in S$, there is a \mbox{transformation} $T$ in $\cols{T}$ such that $T(x) = y$.

\begin{theorem}\label{thm: characterization of homogenuity}
Let $(S, \cols{T})$ be a geometry. Then $\cols{F} = \cset{\set{s}}{s\in S}$ is minimally invariant if and only if $(S, \cols{T})$ is homogeneous.
\end{theorem}
\begin{proof}
The proof of the theorem is straightforward.
\end{proof}

The notion of homogeneity may be restated in terms of group actions. In fact, if $(S, \cols{T})$ is a geometry, then the transformation group $\cols{T}$ acts on the set $S$ of points by evaluation: $T\cdot s = T(s)$ for all $T\in\cols{T}, s\in S$. Thus, $(S, \cols{T})$ is homogeneous if and only if this group action is transitive. From this point of view, we can translate the notion of $n$-transitivity of a group action to the case of geometries. 

\begin{definition}
Let $n\in\N$. A geometry $(S, \cols{T})$ is {\it $n$-transitive} if $S$ has at least $n$ distinct elements and if for each pair of $n$-tuples $(x_1, x_2,\ldots, x_n), (y_1, y_2,\ldots, y_n)\in S^n$ with $x_i\ne x_j$ and $y_i\ne y_j$, $1 \leq i \ne j\leq n$, there exists a transformation $T\in\cols{T}$ such that $T(x_i) = y_i$ for all $i = 1,2,\ldots, n$. 
\end{definition}

In the literature, being $1$-transitive is simply called being transitive and being $2$-transitive is sometimes called being doubly transitive. By definition, we obtain the \mbox{following} proposition immediately.

\begin{proposition}\label{prop: n-transitive --> m-transitive}
If a geometry $(S, \cols{T})$ is $n$-transitive, then it is $m$-transitive for all $m$ with $m\leq n$.
\end{proposition}

The importance of being $n$-transitive lies in the following theorem, which provides one way to construct minimally invariant sets.

\begin{theorem}\label{thm: n-transitive --> minimally invariant}
Let $(S, \cols{T})$ be a geometry with $\abs{S}\geq n$. If $(S, \cols{T})$ is $n$-transitive, then $$\cols{F}_n = \cset{F}{F\subseteq S\textrm{ and }\abs{F} = n}$$ is minimally invariant.
\end{theorem}

The converse to Theorem \ref{thm: n-transitive --> minimally invariant} is true when $n = 1$ by Theorem \ref{thm: characterization of homogenuity}, but not in general true when $n\geq 2$. For instance, let $S = \set{1, 2, 3}$ and let $\cols{T}$ be the alternating group $A_3 = \set{I, \sigma, \tau}$, where $I$ is the identity function on $S$, $\sigma = (1\gap 2\gap 3)$, and $\tau = (1\gap 3\gap 2)$. Note that $\cols{F}_2 = \set{\set{1, 2}, \set{1, 3}, \set{2, 3}}$ is minimally invariant since
\begin{align*}
\set{1, 2} &\overset{I}{\rightarrow} \set{1, 2},\quad \set{1, 2} \overset{\tau}{\rightarrow} \set{1, 3}, \quad \set{1, 2} \overset{\sigma}{\rightarrow} \set{2, 3},\\
\set{1, 3} &\overset{I}{\rightarrow} \set{1, 3},\quad \set{1, 3} \overset{\tau}{\rightarrow} \set{2, 3}, \quad \set{1, 3} \overset{\sigma}{\rightarrow} \set{1, 2},\\
\set{2, 3} &\overset{I}{\rightarrow} \set{2, 3},\quad \set{2, 3} \overset{\tau}{\rightarrow} \set{1, 2}, \quad \set{2, 3} \overset{\sigma}{\rightarrow} \set{1, 3}.
\end{align*}
However, $(S, A_3)$ is not $2$-transitive since there is no transformation in $A_3$ that sends $1$ to $1$ and $2$ to $3$.

We close this section with the definition of being sharply $n$-transitive.
\begin{definition}
A geometry $(S, \cols{T})$ is {\it sharply $n$-transitive} if $S$ has at least $n$ distinct elements and if for each pair of $n$-tuples $(x_1, x_2,\ldots, x_n), (y_1, y_2,\ldots, y_n)\in S^n$ with $x_i\ne x_j$ and $y_i\ne y_j$, $1 \leq i \ne j\leq n$, there exists a unique transformation $T\in\cols{T}$ such that $T(x_i) = y_i$ for all $i = 1,2,\ldots, n$. 
\end{definition}

It is clear that  being sharply $n$-transitive implies being $n$-transitive. However, the reverse implication is not, in general, true. The following theorem gives a sufficient condition, related to fixed points, for the reverse implication to be true. Let $(S, \cols{T})$ be a geometry. For each $T\in\cols{T}$, define
\begin{equation}
\Fix{(T)} = \cset{x\in S}{T(x) = x}.
\end{equation}
Throughout the remainder of this article, $I$ denotes the identity transformation on the set under consideration.

\begin{theorem}
Let $(S, \cols{T})$ be a geometry and let $n\in\N$. If $\abs{\Fix{(T)}}\leq n-1$ for all nonidentity transformations $T\in\cols{T}$, then $(S,\cols{T})$ is $n$-transitive if and only if $(S,\cols{T})$ is sharply $n$-transitive.
\end{theorem}
\begin{proof}
Suppose that $(S,\cols{T})$ is $n$-transitive. Let $(x_1, x_2,\ldots, x_n), (y_1, y_2,\ldots, y_n)\in S^n$ with $x_i\ne x_j$ and $y_i\ne y_j$, $1 \leq i \ne j\leq n$. By definition, there exists a transformation $T$ such that $T(x_i) = y_i$ for all $i = 1,2,\ldots, n$. Suppose that $R\in\cols{T}$ and that $R(x_i) = y_i$ for all $i = 1,2,\ldots, n$. Note that $R^{-1}\circ T\in\cols{T}$ and that $(R^{-1}\circ T)(x_i) = x_i$ for all $i = 1,2,\ldots, n$. Hence, $\abs{\Fix{(R^{-1}\circ T)}} \geq n > n-1$. It follows from the assumption that $R^{-1}\circ T = I$ and so $R = T$. This proves that $(S,\cols{T})$ is sharply $n$-transitive.
\end{proof}

\subsection{Basic knowledge of gyrogroup theory}
In this section, we summarize basic definitions and properties of gyrogroups for easy reference. See \cite{AU2008AHG, TS2016TAG} for more details.

In the case when $\oplus$ is a binary operation on a nonempty set $G$, let $\aut{G}$ be the group of automorphisms of $(G, \oplus)$.

\begin{definition}\label{def: gyrogroup}
A nonempty set $G$, together with a binary operation $\oplus$ on $G$, is called a {\it gyrogroup} if it satisfies the following axioms:
\begin{enumerate}
    \item[(G1)] There exists an element $e\in G$ such that $e\oplus a =
    a$ for all $a\in G$.
    \item[(G2)] For each $a\in G$, there exists an element $b\in G$ such that
$b\oplus a = e$.
    \item[(G3)] For all $a$, $b\in G$, there is an automorphism
$\gyr{a,b}{}\in\aut{G}$ such that
    \begin{equation}\tag{left gyroassociative law} a\oplus (b\oplus c) = (a\oplus b)\oplus\gyr{a,
    b}{c}\end{equation}
    for all $c\in G$.
    \item[(G4)] For all $a$, $b\in G$, $\gyr{a\oplus b, b}{} = \gyr{a,
    b}{}$.\hfill(left loop property)
\end{enumerate}
\end{definition}

The element $e$ in (G1) is indeed a unique two-sided identity of $G$, called the {\it identity element}. The element $b$ in (G2) is indeed a unique two-sided inverse of $a$ in $G$, denoted by $\ominus a$. The automorphism $\gyr{a, b}{}$ in (G3) is called the {\it gyroautomorphism} generated by $a$ and $b$. Let $G$ be a gyrogroup. For each $a\in G$, the function $L_a$ defined by $L_a(x) = a\oplus x$ for all $x\in G$ is a bijection from $G$ to itself, called the {\it left gyrotranslation} by $a$. Left gyrotranslations and gyroautomorphisms are intertwined, as reflected in the following equation:
\begin{equation}
L_a\circ L_b = L_{a\oplus b}\circ \gyr{a, b}{}
\end{equation}
for all $a, b\in G$ (cf. Theorem 10 of \cite{MR3362496}). In any gyrogroup, the {\it left cancellation law} holds: $\ominus a\oplus (a\oplus b) = b$. By the left cancellation law, $L_{\ominus a} = L_a^{-1}$ for all $a\in G$, where $L_a^{-1}$ is the inverse function of $L_a$ with respect to composition. It is clear that $L_e$ is the identity permutation of $G$.

Let $G$ be a gyrogroup. A nonempty subset $H$ of $G$ is called a {\it subgyrogroup} if it is a gyrogroup under the operation of $G$ and $\gyr{a, b}{(H)} = H$ for all $a, b\in H$. If $H$ is a subgyrogroup of $G$ such that $\gyr{a, h}{(H)} = H$ for all $a\in G, h\in H$, then it is called an {\it L-subgyrogroup}. A function $\vphi\colon G\to H$, where $G$ and $H$ are gyrogroups, is a {\it gyrogroup homomorphism} if $\vphi(a\oplus b) = \vphi(a)\oplus\vphi(b)$ for all $a, b\in G$. If $\vphi$ is a gyrogroup homomorphism with domain $G$, then the {\it kernel} of $\vphi$ is defined as $\ker{\vphi} = \cset{a\in G}{\vphi(a) = e}$.

Imposing a length function on a gyrogroup allows us to do analysis on that gyrogroup. Recall the following definition of a normed gyrogroup (cf. Definition 2 of \cite{TS2019MSN}).

\begin{definition}
Let $G$ be a gyrogroup. A function $\norm{\cdot}\colon G\to \R$ is called a {\it gyronorm} on $G$ if the following properties hold:
\begin{enumerate}
\item\label{item: positivity} $\norm{x}\geq 0$  for all $x\in G$ and $\norm{x} = 0$ if and only if $x = e$;\hfill(positivity)
\item\label{item: invariant under taking inverses} $\norm{\ominus x} = \norm{x}$ for all $x\in G$;\hfill(invariant under taking inverses)
\item\label{item: subadditivity} $\norm{x\oplus y}\leq \norm{x}+\norm{y}$ for all $x, y\in G$;\hfill(subadditivity)
\item\label{item: invariant under gyrations}$\norm{\gyr{a, b}{x}} = \norm{x}$ for all $a, b, x\in G$.\hfill(invariant under gyrations)
\end{enumerate}
In this case, $(G, \norm{\cdot})$ is called a {\it normed gyrogroup}.
\end{definition}

Associated to each normed gyrogroup $(G, \norm{\cdot})$ is the metric $d_{\eta}$ defined by
\begin{equation}\label{eqn: gyronorm metric}
d_{\eta}(x, y) = \norm{\ominus x\oplus y}
\end{equation}
for all $x, y\in G$, called the {\it gyronorm metric}. Hence, $(G, d_{\eta})$ forms a metric space. One of the virtues of the gyronorm metric is that it is invariant under left gyrotranslations; that is, $d_{\eta}(a\oplus x, a\oplus y) = d_{\eta}(x, y)$ for all gyrogroup elements $a, x, y$. In other words, the left gyrotranslations of a normed gyrogroup $G$ are isometries of $(G, d_{\eta})$. Furthermore, the gyroautomorphisms of $G$ are isometries of $(G, d_{\eta})$ because they preserve both the gyrogroup operation and the gyronorm. See \cite{TS2019MSN} for more details.

\section{Geometry of gyrogroups via the Erlangen Program}\label{sec: Geometry of gyrogroups}
The study of geometry of abstract gyrogroups in this section is inspired by our three research articles \cite{MR4038385, Suksumran2020505, TRTS2021IRM}. 

In any gyrogroup $G$, there are two fundamental classes of transformations. The first class consists of all left gyrotranslations $L_a\colon x\mapsto  a\oplus x$. The second class consists of all gyroautomorphisms of $G$. None of these forms a group under composition. It is evident that algebraic and geometric properties of a gyrogroup do depend on its left gyrotranslations and its gyroautomorphisms. Therefore, it seems reasonable to study a geometry of an arbitrary gyrogroup using a group of transformations that contains all left gyrotranslations and gyroautomorphisms.

We begin with a parametrization of the symmetric group $\sym{G}$ by left gyro-translations and permutations in $\stab{e}$, where $G$ is a gyrogroup and
\begin{equation}
\stab{e} = \cset{\rho\in\sym{G}}{\rho(e) = e}.
\end{equation}

\begin{theorem}
If $G$ is a gyrogroup, then the symmetric group $\sym{G}$ can be para-meterized by the left gyrotranslations and the permutations of $G$ preserving the gyrogroup identity in the following sense:
\begin{enumerate}
\item\label{item: unique factorization} Each permutation $\sigma$ in $\sym{G}$ has a unique factorization $\sigma = L_a\circ \rho$, where $a\in G$ and $\rho\in \stab{e}$. 
\item For all $a, b\in G, \alpha, \beta\in\stab{e}$, $\gyr{a, \alpha(b)}{}\circ L_{\ominus \alpha(b)}\circ\alpha\circ L_b\circ \beta$ lies in $\stab{e}$ and 
\begin{equation}\label{eqn: composition law in Sym(G)}
(L_a\circ \alpha)\circ(L_b\circ\beta) = L_{a\oplus\alpha(b)}\circ (\gyr{a, \alpha(b)}{}\circ L_{\ominus \alpha(b)}\circ\alpha\circ L_b\circ \beta).
\end{equation}
\end{enumerate}
\end{theorem}
\begin{proof}
Part \eqref{item: unique factorization} was proved in Theorem 11 of \cite{MR3362496}. Let $a, b\in G, \alpha, \beta\in\stab{e}$. Note that 
$$
(L_{\ominus \alpha(b)}\circ\alpha\circ L_b\circ \beta)(e) = \ominus\alpha(b)\oplus(\alpha(b\oplus\beta(e))) = e.
$$
Hence, $(\gyr{a, \alpha(b)}{}\circ L_{\ominus \alpha(b)}\circ\alpha\circ L_b\circ \beta)(e) = \gyr{a, \alpha(b)}{e} = e$. Note that $$\alpha\circ L_b = L_{\alpha(b)}\circ (L_{\ominus\alpha(b)}\circ \alpha\circ L_b)$$ since $L_{\ominus\alpha(b)} = L^{-1}_{\alpha(b)}$. By Equation (12) of \cite{MR3362496}, $L_a\circ L_{\alpha(b)} = L_{a\oplus\alpha(b)}\circ\gyr{a, \alpha(b)}{}$. Therefore, $(L_a\circ \alpha)\circ(L_b\circ\beta) = L_{a\oplus\alpha(b)}\circ (\gyr{a, \alpha(b)}{}\circ L_{\ominus \alpha(b)}\circ\alpha\circ L_b\circ \beta)$ as desired.
\end{proof}

Note that Equation \eqref{eqn: composition law in Sym(G)} generalizes the notion of a gyrosemidirect product. In fact, if $\alpha$ is an automorphism of $G$, then the commutation relation (cf. Equation (14) of \cite{MR3362496}) implies that $\alpha\circ L_b = L_{\alpha(b)}\circ \alpha$ and so Equation \eqref{eqn: composition law in Sym(G)} reduces to the usual composition law of the gyrosemidirect product group $G\rtimes_{\rm gyr} \aut{G}$ (cf. Section 2.6 of \cite{AU2008AHG}).

Let $G$ be a gyrogroup. According to the fact that $\gyr{a, b}{}$ is an automorphism of $G$ for all $a, b\in G$, we can let $\GYR{G}$ be the subgroup of $\sym{G}$ generated by all the gyroautomorphisms of $G$; that is,
\begin{equation}
\GYR{G} = \gen{\gyr{a, b}{}\colon a, b\in G}.
\end{equation}
Clearly, $\GYR{G}$ is a subgroup of $\aut{G}$. Furthermore, $G$ forms a group under the gyrogroup operation if and only if $\GYR{G} = \set{I}$. Note that, as $\igyr{a,b}{} = \gyr{b, a}{}$ for all $a, b\in G$, we obtain
\begin{equation}
\GYR{G} = \cset{\gyr{a_1, b_1}{}\circ \gyr{a_2, b_2}{}\circ\cdots\circ \gyr{a_n, b_n}{}}{a_i, b_i\in G, i = 1,2,\ldots, n}.
\end{equation}

Let $G$ be a gyrogroup. Define 
\begin{eqnarray}
\hat{G} &=& \cset{L_a}{a\in G},\\
\Gamma_M &=& \cset{L_a\circ \tau}{a\in G, \tau\in\aut{G}},\\ 
\Gamma_m &=&  \cset{L_a\circ \gamma}{a\in G, \gamma\in\GYR{G}}. \label{eqn: Gamma m}
\end{eqnarray}
Using the fact that
\begin{eqnarray*}
(L_a\circ \alpha)\circ(L_b\circ \beta) &=& L_{a\oplus \alpha(b)}\circ (\gyr{a, \alpha(b)}{}\circ\alpha\circ\beta)\\
(L_a\circ \alpha)^{-1} &=& L_{\alpha^{-1}(\ominus a)}\circ\alpha^{-1}
\end{eqnarray*}
for all $a, b\in G, \alpha, \beta\in\aut{G}$, one can verify that $\Gamma_M$ and $\Gamma_m$ form subgroups of $\sym{G}$ and hence are groups under composition. Indeed, the group $\Gamma_M$ is isomorphic to the {\it \mbox{gyroholomorph}} of $G$ (cf. Definition 2.29 of \cite{AU2008AHG}). As evidenced in the literature, \mbox{being} invariant under translations is a nice mathematical property. Therefore, we will focus on geometries whose groups of transformations contain left gyrotranslations. The next proposition shows that $\Gamma_m$ is a minimal subgroup of $\sym{G}$ that contains all left gyrotranslations.

\begin{proposition}
Let $G$ be a gyrogroup and let $\Sigma$ be a subgroup of $\sym{G}$. Then $\hat{G}\subseteq \Sigma$ if and only if $\Gamma_m\subseteq \Sigma$.
\end{proposition}
\begin{proof}
Suppose that $\hat{G}\subseteq \Sigma$. Let $a, b\in G$. By Equation (12) of \cite{MR3362496}, $$\gyr{a, b}{} = L^{-1}_{a\oplus b}\circ L_a\circ L_b$$ and so $\gyr{a, b}{}\in \Sigma$. This implies that $\GYR{G}\subseteq \Sigma$. It follows that $\Gamma_m\subseteq \Sigma$ since $\hat{G}\GYR{G}\subseteq \Sigma$ by the closure property. The converse is trivial since $\hat{G}\subseteq \Gamma_m$.
\end{proof}

Throughout the remainder of this section, we will explore geometries of gyrogroups with respect to the groups $\Gamma_m$ and $\Gamma_M$.

\subsection{The property of $n$-transitivity}
In this section, we focus on the property of transitivity and $n$-transitivity of gyro-groups. These results are abstract versions of Corollary 6 in \cite{MR3646419}, Theorem 2.7 in \cite{MR4038385}, Theorem 2.6 in \cite{Suksumran2020505}, and Theorem 1 in \cite{TRTS2021IRM}.

\begin{theorem}\label{thm: homogenuity of Gamma m}
Let $G$ be a gyrogroup. The geometry $(G, \Gamma_m)$ is homogeneous (and hence transitive): if $x$ and $y$ are points of $G$, then there exists a transformation $T$ in $\Gamma_m$ such that $T(x) = y$.
\end{theorem}
\begin{proof}
Let $x, y\in G$. Define $T= L_y\circ L_{\ominus x}$. By Theorem 18 (3) of \cite{TS2016TAG}, 
$$
T = L_{y\ominus x}\circ\gyr{y, \ominus x}{}
$$
and so $T\in \Gamma_m$. Furthermore,  $T(x) = (L_y\circ L_{\ominus x})(x) =y\oplus (\ominus x\oplus x) = y\oplus e = y$.
\end{proof}

\begin{corollary}\label{cor: homogenuity of Gamma M}
Let $G$ be a gyrogroup. The geometry $(G, \Gamma_M)$ is homogeneous (and hence transitive).
\end{corollary}
\begin{proof}
The corollary follows from the fact that $\Gamma_m\subseteq \Gamma_M$.
\end{proof}

\begin{corollary}
If $G$ is a gyrogroup, then $\cset{\set{a}}{a\in G}$ is a minimally invariant set of $(G, \Gamma_m)$ and $(G, \Gamma_M)$.
\end{corollary}
\begin{proof}
The corollary follows from Theorems \ref{thm: characterization of homogenuity} and \ref{thm: homogenuity of Gamma m} and Corollary \ref{cor: homogenuity of Gamma M}.
\end{proof}

In view of Theorem \ref{thm: homogenuity of Gamma m}, it is natural to ask whether the geometry $(G, \Gamma_m)$ is sharply transitive. The answer is \qt{no} in the case when the gyrogroup possesses a nontrivial gyro-automorphism and the answer is \qt{yes} otherwise. By a {\it nondegenerate gyrogroup} we mean a gyrogroup with a nonidentity gyroautomorphism. Hence, every \mbox{nondegenerate} gyrogroup is a gyrogroup that is not a group under the gyrogroup operation. The following theorem ensures the existence of a nonidentity transformation in $\Gamma_m$ that leaves a given point fixed if $G$ is a nondegenerate gyrogroup. 

\begin{theorem}\label{thm: isotropic geometry, gyrogroup with length}
If $G$ is a nondegenerate gyrogroup, then for each point $p\in G$, there exists a nonidentity transformation $T$ in $\Gamma_m$ such that $T(p) = p$.
\end{theorem}
\begin{proof}
Let $p$ be an arbitrary point of $G$. Let $\gamma$ be a {\it nonidentity} gyroautomorphism of $G$. Then $\gamma(e) = e$. Define $T = L_p\circ \gamma\circ L_{\ominus p}$. Note that $T(p) = p$. By Equation (55) of \cite{TS2016TAG} and Theorem 18 (3) of \cite{TS2016TAG},
$$
T =  L_p\circ \gamma\circ L_{\ominus p} = L_p\circ(L_{\gamma(\ominus p)}\circ \gamma) = L_{p\,\oplus\,\gamma(\ominus p)}\circ (\gyr{p, \gamma(\ominus p)}{}\circ \gamma)
$$
and so $T\in \Gamma_m$. Note that $T$ is not the identity transformation of $G$; otherwise, we would have $I = L_p\circ\gamma\circ L_{\ominus p} = L_p\circ\gamma\circ L_p^{-1}$ and would have $\gamma = I$, a contradiction.
\end{proof}

\begin{theorem}\label{thm: nondegenerate gyrogroup transitive but not sharply}
If $G$ is a nondegenerate gyrogroup, then $(G, \Gamma_m)$ is transitive but not sharply transitive.
\end{theorem}
\begin{proof}
By Theorem \ref{thm: homogenuity of Gamma m}, $(G, \Gamma_m)$ is transitive. Let $x, y\in G$ and let $T$ be as in the proof of Theorem \ref{thm: homogenuity of Gamma m}. Then $T(x) = y$.  By Theorem \ref{thm: isotropic geometry, gyrogroup with length}, there is a nonidentity transformation $R$ in $\Gamma_m$ such that $R(x) = x$. Note that $T\circ R\ne T$ since otherwise $R$ would be the identity transformation. Note that $(T\circ R)(x) = T(R(x)) = T(x) = y$. Therefore, $(G, \Gamma_m)$ is not sharply transitive.
\end{proof}

\begin{corollary}
If $G$ is a nondegenerate gyrogroup, then $(G, \Gamma_M)$ is transitive but not sharply transitive.
\end{corollary}
\begin{proof}
The corollary follows from the fact that $\Gamma_m\subseteq \Gamma_M$.
\end{proof}

In the case when the gyroautomorphisms of $G$ are trivial, $G$ becomes a group and $\Gamma_m$ reduces to the group $\cols{T}$ of left translations. In this case, the translational geometry $(G, \cols{T})$ is sharply transitive.

\begin{theorem}
If $G$ is a group, then the translational geometry $(G, \cols{T})$ is sharply transitive.
\end{theorem}
\begin{proof}
Let $x, y\in G$. Suppose that $a, b\in G$ and that $L_a(x) = y$ and $L_b(x) = y$. Then $a\oplus x = y = b\oplus x$. By the right cancellation law in groups, $a = b$ and so $L_a = L_b$.
\end{proof}

Next, we consider the question of whether the geometry $(G, \Gamma_m)$ is $n$-transitive. Let $G$ be a gyrogroup. Set $\check{G}= G\setminus\set{e}$, where $e$ is the gyrogroup identity, and let
\begin{equation}
\GYR{\check{G}} = \cset{\res{\gamma}{\check{G}}}{\gamma\in\GYR{G}}.
\end{equation}
Here, $\res{\gamma}{\check{G}}$ is the restriction of $\gamma$ to $\check{G}$. Since every transformation in $\GYR{G}$ is bijective and preserves $e$, it follows that $\GYR{\check{G}}$ forms a group under composition. Therefore, we can speak of the geometry $(\check{G}, \GYR{\check{G}})$. Intransitive geometry of $(\check{G}, \GYR{\check{G}})$ leads to $n$-intransitive geometry of $(G,\Gamma_m)$ for all $n\geq 2$, as shown in the following theorem.

\begin{theorem}\label{thm: sufficient condition to be not n-transitive}
Let $G$ be a nontrivial gyrogroup. If $(\check{G}, \GYR{\check{G}})$ is not transitive, then $(G, \Gamma_m)$ is not $n$-transitive for all $n\geq 2$. 
\end{theorem}
\begin{proof}
We prove the contrapositive. Suppose that $(G, \Gamma_m)$ is $n$-transitive for some $n\geq 2$. By Proposition \ref{prop: n-transitive --> m-transitive}, $(G, \Gamma_m)$ is $2$-transitive. Let $x, y\in\check{G}$. Then $x\ne e$ and $y\ne e$. By definition, there is a transformation $T = L_a\circ \gamma$ in $\Gamma_m$, where $a\in G$ and $\gamma\in\GYR{G}$, such that $T(e) = e$ and $T(x) = y$. Note that $e = T(e) = a\oplus \gamma(e) = a\oplus e = a$ and so $L_a = I$. Hence, $T = \gamma\in\GYR{G}$. It follows that $\res{T}{\check{G}}\in\GYR{\check{G}}$. Furthermore, $\res{T}{\check{G}}(x) = y$. This shows that $(\check{G}, \GYR{\check{G}})$ is transitive.
\end{proof} 

As an application of Theorem \ref{thm: sufficient condition to be not n-transitive}, we give a sufficient condition for the geometry $(G,\Gamma_m)$ to be $n$-intransitive for all $n\geq 2$. Let $G$ be a gyrogroup. Recall that the {\it right nucleus} of $G$, denoted by $N_r(G)$, is defined by 
\begin{equation}
N_r(G) = \cset{c\in G}{\gyr{a, b}{c} = c\textrm{ for all }a, b\in G}.
\end{equation}
It is a special subgyrogroup of $G$ that forms a group under the gyrogroup operation; see Section 3 of \cite{TS2016SSG} for more details.

\begin{theorem}\label{thm: right nucleus nontrivial --> not n-transitive}
Let $G$ be a gyrogroup with at least three distinct elements. If $N_r(G)$ is nontrivial, then $(G, \Gamma_m)$ is not $n$-transitive for all $n\geq 2$. 
\end{theorem}
\begin{proof}
Suppose that $N_r(G)$ is nontrivial. In light of Theorem \ref{thm: sufficient condition to be not n-transitive}, it suffices to show that $(\check{G}, \GYR{\check{G}})$ is not transitive. By assumption, there is an element $x\in N_r(G)$ with $x\ne e$. Since $\abs{G}\geq 3$, we can chose an element $y\in G\setminus\set{e, x}$. Then $y\in\check{G}$. Since $x\in N_r(G)$, it follows that $\gamma(x) = x$ for all $\gamma\in\GYR{G}$. This implies that there is no transformation $T$ in $\GYR{\check{G}}$ such that $T(x) = y$ since otherwise we would have $x = y$.
\end{proof}

Theorem \ref{thm: right nucleus nontrivial --> not n-transitive} provides a convenient way to confirm that a given gyrogroup is not $n$-transitive for all $n\geq 2$. Some concrete examples are exhibited below.

\begin{example}
In the gyrogroup $K_{16} =  \set{0,1,2,3,4,5,6,7,8,9,10,11,12,13,14,15}$ (cf. \cite[p. 41]{AU2002BEA}), the only nontrivial gyroautomorphism of $K_{16}$ is given in terms of cycle notation by
\begin{equation}
A = (8\gap 9)(10\gap 11)(12\gap 13)(14\gap 15).
\end{equation}
Hence, $N_r(K_{16}) = \set{0,1,2,3,4,5,6,7}$. By Theorem \ref{thm: right nucleus nontrivial --> not n-transitive}, $(K_{16}, \Gamma_m)$ is not $n$-transitive for all $n\geq 2$. 
\end{example}

\begin{example}
In the gyrogroup $G_{8} =\set{0,1,2,3,4,5,6,7}$ (cf. \cite[p. 404]{TS2016TAG}), the only nontrivial gyroautomorphism of $G_{8}$ is given in terms of cycle notation by
\begin{equation}
A = (4\gap 6)(5\gap 7).
\end{equation}
Hence, $N_r(G_{8}) = \set{0,1,2,3}$. By Theorem \ref{thm: right nucleus nontrivial --> not n-transitive}, $(G_{8}, \Gamma_m)$ is not $n$-transitive for all $n\geq 2$. 
\end{example}

\begin{example}
In the gyrogroup $G_{15} =\set{0,1,2,3,4,5,6,7,8,9,10,11,12,13,14}$  (cf. \cite[p. 432]{TS2016TAG}), the nontrivial gyroautomorphisms of $G_{15}$ are given in terms of cycle notation by
\begin{eqnarray}\label{eqn: gyration of G15}
\begin{split}
A &= (1\gap 7 \gap 5 \gap 10 \gap 6)(2 \gap 3 \gap 8
\gap 11 \gap 14),\\
B &= (1\gap 6 \gap 10 \gap 5 \gap 7)(2 \gap 14 \gap 11
\gap 8 \gap 3),\\
C &= (1\gap 10 \gap 7 \gap 6 \gap 5)(2 \gap 11 \gap 3
\gap 14 \gap 8),\\
D &= (1\gap 5 \gap 6 \gap 7 \gap 10)(2 \gap 8 \gap 14
\gap 3 \gap 11).
\end{split}
\end{eqnarray}
Hence, $N_r(G_{15}) = \set{0, 4, 9, 12, 13}$. By Theorem \ref{thm: right nucleus nontrivial --> not n-transitive}, $(G_{15}, \Gamma_m)$ is not $n$-transitive for all $n\geq 2$. 
\end{example}

The converse to Theorem \ref{thm: right nucleus nontrivial --> not n-transitive} is not true, as shown in the next example. Recall that the {\it complex M\"{o}bius gyrogroup} \cite{AU2008FMG} consists of the open unit disk $\D = \cset{z\in\C}{\abs{z}<1}$, together with M\"{o}bius addition $\oplus_M$ defined by
\begin{equation}
w\oplus_M z = \dfrac{w+z}{1+\bar{w}z}
\end{equation}
for all $w, z\in\D$. Its identity is $0$ and  its gyroautomorphisms, which are disk rotations, are described by 
\begin{equation}
\gyr{a, b}{z} = \dfrac{1+a\bar{b}}{1+ \bar{a}b}z
\end{equation}
for all $a, b, z\in\D$.

\begin{example}\label{exm: complex Mobius gyrogroup}
The geometry $(\D, \Gamma_m)$ is not $n$-transitive for all $n\geq 2$. Furthermore, $N_r(\D) = \set{0}$. In fact, we chose $x, y\in\D\setminus\set{0}$ with $\abs{x}\ne \abs{y}$; for example, $x = 0.5$ and $y = 0.7$. Suppose that $\gamma = \gyr{a_1, b_1}{}\circ \gyr{a_2, b_2}{}\circ \cdots \circ \gyr{a_k, b_k}{}$. Then $\gamma(x)\ne y$ since otherwise we would have $\abs{y} = \abs{\gamma(x)} = \abs{x}$, a contradiction. This shows that $(\check{\D}, \GYR{\check{\D}})$ is not transitive. By Theorem \ref{thm: sufficient condition to be not n-transitive}, $(\D, \Gamma_m)$ is not $n$-transitive for all $n\geq 2$. Note that if $\gyr{a, b}{}\ne I$, then $\dfrac{1+a\bar{b}}{1+\bar{a}b}\ne 1$. Hence, $\gyr{a, b}{c} = c$ implies $\dfrac{1+a\bar{b}}{1+\bar{a}b}c = c$, which in turn implies $c = 0$. Therefore, $N_r(\D) = \set{0}$.
\end{example}

So far we have no example of a gyrogroup $G$ such that $(G, \Gamma_m)$ is $n$-transitive for some $n\geq 2$. This leads to the following question.

\begin{question}
Is there a (nondegenerate) gyrogroup $G$ such that $(G, \Gamma_m)$ is $n$-transitive for some $n\geq 2$?
\end{question}

We emphasize that any gyrogroup that might be an example of the previous question must have the trivial right nucleus by Theorem \ref{thm: right nucleus nontrivial --> not n-transitive}.

\subsection{Coset spaces as invariant sets}
Recall that if $\Gamma$ is a group, written multiplicatively, and if $\Xi$ is a subgroup of $\Gamma$, then the left multiplication function $\lambda_a$ defined by $\lambda_a(g) = ag$ for all $g\in\Gamma$ sends the left coset $g\Xi$ to the left coset $(ag)\Xi$.  Geometrically, the coset space $\Gamma/\Xi = \cset{g\Xi}{g\in \Gamma}$ is an invariant set of the geometry $(\Gamma, \cols{T})$, where $\cols{T} = \cset{\lambda_a}{a\in \Gamma}$.  This is not the case for gyrogroups because of the presence of gyroautomorphisms. In fact, there is an example of a gyrogroup $G$ such that its coset space is not an invariant set of $(G, \Gamma_m)$. Therefore, we seek to find some necessary and sufficient conditions for the left coset space of a gyrogroup to be an invariant set. This leads to the following definition.

\begin{definition}
A subgyrogroup $H$ of a gyrogroup $G$ is {\it strong}, denoted by $H\leq_S G$, if $\gyr{a, b}{(H)} = H$ for all $a, b\in G$. It is called a {\it characteristic subgyrogroup} of $G$ if $\tau(H) = H$ for all $\tau\in\aut{G}$.
\end{definition}

It is clear that every characteristic subgyrogroup is strong. By Proposition 32 (3) of \cite{TS2016TAG}, every normal subgyrogroup is strong. Moreover, by definition, every strong subgyrogroup is an L-subgyrogroup. The following theorem gives a necessary and sufficient condition for the coset space to be an invariant set of the geometry $(G, \Gamma_m)$, where $G$ is an arbitrary gyrogroup.

\begin{theorem}\label{thm: characterization of coset to be invariant}
Let $H$ be a subgyrogroup of a gyrogroup $G$. Then the collection of left cosets $G/H = \cset{a\oplus H}{a\in G}$ is an invariant set of the geometry $(G, \Gamma_m)$ if and only if $H$ is  strong in $G$.
\end{theorem}
\begin{proof}
($\If$) Let $a, b\in G$. By definition, $\gyr{a, b}{(H)}\in G/H$ and so $\gyr{a, b}{(H)} = c\oplus H$ for some $c\in G$. Since $e = \gyr{a, b}{e}\in \gyr{a, b}{(H)}$, it follows that $e = c\oplus h$ for some $h\in H$. Thus, $c = \ominus h$ and so $c\oplus H = \ominus h\oplus H = H$ by the closure property. This proves that $H\leq_S G$.

($\Fi$) Let $g\in G$. Using the left gyroassociative law, we obtain that
$$
a\oplus (g\oplus h) = (a\oplus g)\oplus\gyr{a, g}{h}
$$
for all $a\in G, h\in H$. It follows that $$L_a(g\oplus H) = a\oplus(g\oplus H) = (a\oplus g)\oplus\gyr{a, g}{(H)} = (a\oplus g)\oplus H$$ for all $a\in G$. Hence, $L_a(g\oplus H)\in G/H$ for all $a\in G$. Let $a, b\in G$. Since $$\gyr{a, b}{(g\oplus h)} = (\gyr{a, b}{g})\oplus(\gyr{a, b}{h})$$ for all $h\in H$, it follows that $$\gyr{a, b}{(g\oplus H)} = (\gyr{a, b}{g})\oplus \gyr{a, b}{(H)} = (\gyr{a, b}{g})\oplus H.$$ Hence, $\gyr{a, b}{(g\oplus H)}\in G/H$ for all $a, b\in G$. Let $T\in\Gamma_m$. Then $$T = L_a\circ\gyr{a_1, b_1}{}\circ \gyr{a_2, b_2}{}\circ\cdots\circ \gyr{a_n, b_n}{}$$ for some $a, a_i, b_i\in G$. It follows that $T(g\oplus H)\in G/H$, which completes the proof.
\end{proof}

Similarly, we obtain a necessary and sufficient condition for the coset space to be an invariant set of the geometry $(G, \Gamma_M)$.

\begin{theorem}\label{thm: characterization of coset to be invariant, Gamma M}
Let $H$ be a subgyrogroup of a gyrogroup $G$. Then the collection of left cosets $G/H$ is an invariant set of the geometry $(G, \Gamma_M)$ if and only if $H$ is  characteristic in $G$.
\end{theorem}
\begin{proof}
The proof of this theorem is similar to that of Theorem \ref{thm: characterization of coset to be invariant}.
\end{proof}

The following two corollaries show how to construct some minimally invariant sets of the geometries $(G, \Gamma_m)$ and $(G, \Gamma_M)$, respectively.

\begin{corollary}\label{cor: G/H is  a minimally invariant set when H is strong}
Let $H$ be a subgyrogroup of a gyrogroup $G$. If $H$ is strong in $G$, then $G/H $ is a minimally invariant set of the geometry $(G, \Gamma_m)$.
\end{corollary}
\begin{proof}
According to Theorem \ref{thm: characterization of minimally invariant set II}, it suffices to show that any two figures in $G/H$ are congruent. Let $X, Y\in G/H$. Then $X = x\oplus H$ and $Y = y\oplus H$ for some $x, y\in G$. Set $T = L_y\circ L_{\ominus x}$. As in the proof of Theorem \ref{thm: homogenuity of Gamma m}, $T\in\Gamma_m$. Since
$$
T(X) = T(x\oplus H) = y\oplus (\ominus x\oplus (x\oplus H)) = y\oplus H = Y,
$$
it follows that $X\cong Y$.
\end{proof}

\begin{corollary}
Let $H$ be a subgyrogroup of a gyrogroup $G$. If $H$ is characteristic in $G$, then $G/H $ is a minimally invariant set of the geometry $(G, \Gamma_M)$.
\end{corollary}
\begin{proof}
The corollary follows from the fact that $\Gamma_m\subseteq \Gamma_M$.
\end{proof}

Due to the fact that the kernel of a gyrogroup homomorphism is invariant under all gyroautomorphisms, we obtain the following theorem immediately.

\begin{theorem}
Let $\vphi$ be a homomorphism of a gyrogroup $G$. Then $\ker{\vphi}$ is a strong subgyrogroup of $G$. Consequently, $G/\ker{\vphi}$ is a minimally invariant set of the geometry $(G, \Gamma_m)$.
\end{theorem}
\begin{proof}
By Proposition 23 of \cite{MR3362496}, $\gyr{a, b}{(\ker{\vphi})} = \ker{\vphi}$ for all $a, b\in G$. Hence, $\ker{\vphi}$ is a strong subgyrogroup of $G$ and so the theorem follows directly from Corollary \ref{cor: G/H is  a minimally invariant set when H is strong}.
\end{proof}

We close this section with an easy characterization of minimally invariant sets. Let $G$ be a gyrogroup. For each nonempty subset $F$ of $G$, define
\begin{equation}
\cols{F}_m(F) = \cset{g\oplus \gamma(F)}{g\in G, \gamma\in\GYR{G}}
\end{equation}
and
\begin{equation}
\cols{F}_M(F) = \cset{g\oplus \tau(F)}{g\in G, \tau\in\aut{G}}.
\end{equation}

\begin{theorem}
Let $G$ be a gyrogroup and let $\cols{F}$ be a nonempty set of figures in $G$.
\begin{enumerate}
\item\label{item: characterization of minimally invariant set, m} $\cols{F}$ is minimally invariant in $(G, \Gamma_m)$ if and only if $\cols{F} = \cols{F}_m(F) $ for some nonempty subset $F$ of $G$.
\item\label{item: characterization of minimally invariant set, M} $\cols{F}$ is minimally invariant in $(G, \Gamma_M)$ if and only if $\cols{F} = \cols{F}_M(F)$ for some nonempty subset $F$ of $G$.
\end{enumerate}
\end{theorem}
\begin{proof}
To prove part \eqref{item: characterization of minimally invariant set, m}, suppose that $\cols{F}$ is minimally invariant in $(G, \Gamma_m)$. Since $\cols{F} \ne\emptyset$, there is a figure $F\in\cols{F}$. As in the proof of Theorem \ref{thm: characterization of minimally invariant set}, $$\cols{F} = [F] = \cset{T(F)}{T\in\Gamma_m} = \cset{a\oplus \gamma(F)}{a\in G, \gamma\in\GYR{G}} = \cols{F}_m(F).$$
Conversely, let $\emptyset \ne F\subseteq G$. Then 
$$
\cols{F}_m(F) = \cset{a\oplus \gamma(F)}{a\in G, \gamma\in\GYR{G}} = \cset{T(F)}{T\in\Gamma_m} = [F]
$$ 
and so $\cols{F}_m(F)$ is minimally invariant. The proof of part \eqref{item: characterization of minimally invariant set, M} can be done in a similar fashion.
\end{proof}

\begin{corollary}
Let $G$ be a gyrogroup and let $F$ be a nonempty subset of $G$.
\begin{enumerate}
\item If $\gyr{a, b}{(F)} = F$ for all $a, b\in G$, then $G/F = \cset{g\oplus F}{g\in G}$ is a minimally invariant set of $(G, \Gamma_m)$.
\item If $\tau(F) = F$ for all $\tau\in\aut{G}$, then $G/F = \cset{g\oplus F}{g\in G}$ is a minimally \mbox{invariant} set of $(G, \Gamma_M)$.
\end{enumerate}
\end{corollary}

\section{Geometry of abstract normed gyrogroups}\label{sec: normed gyrogroups}

As evidenced in the literature, gyrogroup structures have closed relationships with M\"{o}bius transformations. One of the most beautiful properties of M\"{o}bius transformations on the extended complex plane $\overline{\C} = \C\cup\set{\infty}$ is shown in the {\it fundamental \mbox{theorem} of M\"{o}bius transformations}. It states that there is a unique M\"{o}bius transformation \mbox{sending} any three distinct points of $\overline{\C} $ to any three distinct points of $\overline{\C} $. Geometrically, $(\overline{\C}, \cols{M})$, where $\cols{M}$ is the group of M\"{o}bius transformations, is sharply $3$-transitive. The \mbox{notion} of a normed gyrogroup may be viewed as a suitable abstract version of the M\"{o}bius gyro-group. This motivates us to study geometry of abstract normed gyrogroups, following the Erlangen Program.

By Theorem \ref{thm: homogenuity of Gamma m}, every normed gyrogroup is transitive. By Theorem \ref{thm: nondegenerate gyrogroup transitive but not sharply}, \mbox{every} nondegenerate normed gyrogroup is not sharply transitive. In light of \mbox{Example} \ref{exm: complex Mobius gyrogroup}, we show that any normed gyrogroup having two nonidentity elements of different norms is not $n$-transitive for all $n\geq 2$.  

\begin{theorem}\label{thm: two non-equal norm implies not n-transitive}
Let $G$ be a normed gyrogroup. If $G$ contains two nonidentity \mbox{elements} $x$ and $y$ such that $\norm{x}\ne\norm{y}$, then $(G,\Gamma_m)$ is not $n$-transitive for all $n\geq 2$.
\end{theorem}
\begin{proof}
According to Proposition \ref{prop: n-transitive --> m-transitive}, it suffices to show that $(G, \Gamma_m)$ is not $2$-transitive. Assume to the contrary that $(G, \Gamma_m)$ is $2$-transitive. By hypothesis, there are elements $x, y\in G\setminus\set{e}$ with $\norm{x}\ne \norm{y}$. By definition, there is a transformation in $\Gamma_m$ of the form $T = L_a\circ \gamma$, where $a\in G$ and $\gamma = \gyr{a_1, b_1}{}\circ \gyr{a_2, b_2}{}\circ \cdots \circ \gyr{a_k, b_k}{}$, such that $T(e) = e$ and $T(x) = y$. This implies $a =e$. Hence, $L_a =I$ and so $T = \gamma$. Since each $\gyr{a_i, b_i}{}$ preserves the gyronorm, it follows that $\norm{y} = \norm{T(x)} = \norm{\gamma(x)} = \norm{x}$, a contradiction. Therefore, $(G, \Gamma_m)$ is not $2$-transitive.
\end{proof}

In the literature, two prominent examples of concrete normed gyrogroups are (Euclidean) Einstein gyrogroups (cf. Section 3.8 of \cite{AU2008AHG}) and (Euclidean) M\"{o}bius gyrogroups (cf. Section 3.5 of \cite{AU2008AHG}). Let $\B$ denote the open unit ball in $n$-dimensional Euclidean space $\R^n$; that is,
\begin{equation}
\B = \cset{\vec{v}\in\R^n}{\norm{\vec{v}}<1},
\end{equation}
where $\norm{\cdot}$ is the {\it Euclidean} norm on $\R^n$. The {\it Einstein gyrogroup} consists of $\B$, together with Einstein addition $\oplus_E$ defined by
\begin{equation}\label{eqn: Einstein addition}
\vec{u}\oplus_E\vec{v} =
\dfrac{1}{1+\gen{\vec{u},\vec{v}}}\Bp{\vec{u}+
\dfrac{1}{\gamma_{\vec{u}}}\vec{v} +
\dfrac{\gamma_{\vec{u}}}{1+\gamma_{\vec{u}}}\gen{\vec{u},\vec{v}}\vec{u}}
\end{equation}
for all $\vec{u}, \vec{v}\in\B$, where $\gamma_{\vec{u}}$ is the {\it Lorentz factor} given by 
$\gamma_{\vec{u}} = \dfrac{1}{\sqrt{1-\norm{\vec{u}}^2}}$. The Einstein gyrogroup becomes a normed gyrogroup when it is endowed with the gyronorm given by
\begin{equation}\label{eqn: length function on Einstein gyrogroup}
\norm{\vec{v}}_E = \tanh^{-1}{\norm{\vec{v}}},\qquad \vec{v}\in\B,
\end{equation}
where $\tanh^{-1}$ is the inverse of the hyperbolic tangent function on $\R$ (cf. Theorem 2 of \cite{TS2019MSN}). 

\begin{corollary}
The geometry $((\B, \oplus_E), \Gamma_m)$ of the Einstein gyrogroup is not $n$-transitive for all $n\geq 2$. Furthermore, it is transitive but not sharply transitive.
\end{corollary}
\begin{proof}
The first statement follows from Theorem \ref{thm: two non-equal norm implies not n-transitive}. The second statement follows from Theorem \ref{thm: nondegenerate gyrogroup transitive but not sharply}.
\end{proof}

The {\it M\"{o}bius gyrogroup} consists of $\B$, together with M\"{o}bius addition $\oplus_M$ defined by
\begin{equation}\label{eqn: Mobius addition}
\vec{u}\oplus_M\vec{v} = \dfrac{(1 + 2\gen{\vec{u},\vec{v}} +
\norm{\vec{v}}^2)\vec{u} + (1 - \norm{\vec{u}}^2)\vec{v}}{1 +
2\gen{\vec{u},\vec{v}} + \norm{\vec{u}}^2\norm{\vec{v}}^2}
\end{equation}
for all $\vec{u}, \vec{v}\in\B$. The M\"{o}bius gyrogroup becomes a normed gyrogroup when it is endowed with the gyronorm given by
\begin{equation}\label{eqn: length function on Mobius gyrogroup}
\norm{\vec{v}}_M= \dfrac{1}{2}\tanh^{-1}\Bp{\dfrac{2\norm{\vec{v}}}{1+\norm{\vec{v}}^2}},\qquad \vec{v}\in\B
\end{equation}
(cf. Theorem 3 of \cite{TS2019MSN}).

\begin{corollary}
The geometry $((\B, \oplus_M), \Gamma_m)$ of the M\"{o}bius gyrogroup is not $n$-transitive for all $n\geq 2$. Furthermore, it is transitive but not sharply transitive.
\end{corollary}
\begin{proof}
The first statement follows from Theorem \ref{thm: two non-equal norm implies not n-transitive}. The second statement follows from Theorem \ref{thm: nondegenerate gyrogroup transitive but not sharply}.
\end{proof}

In any normed gyrogroup, we can measure the size of a gyrogroup element so that we can define rotation transformations, similar to rotations in Euclidean geometry.

\begin{definition}
Let $G$ be a normed gyrogroup and let $p$ be a point in $G$. A transformation $\rho$ in $\sym{G}$ is called a {\it rotation} of $G$ {\it about} $p$ if $\rho(p) = p$ and $\norm{\rho(x)} = \norm{x}$ for all $x\in G$. 
\end{definition}

Note that all the gyroautomorphisms of a normed gyrogroup $G$ are rotations about the identity element because $\gyr{a, b}{e} = e$ and $\norm{\gyr{a, b}{x}} = \norm{x}$ for all $a, b, x\in G$. This implies that if $\gamma\in\GYR{G}$, then $\gamma$ is a rotation about $e$ because $\gamma$ is a finite composite of rotations about $e$. It is not difficult to check that the set of rotations about a given point forms a group under composition. Throughout the remainder of this section, $d_{\eta}$ is the gyronorm metric described by Equation \eqref{eqn: gyronorm metric}.

\begin{proposition}
Let $G$ be a normed gyrogroup. If $T$ is an isometry of $(G, d_\eta)$ that preserves the identity, then $T$ preserves the gyronorm.
\end{proposition}
\begin{proof}
Direct computation shows that
$$
\norm{x} = d_\eta(e, x) = d_\eta(T(e), T(x)) = d_\eta(e, T(x)) = \norm{T(x)}
$$
for all $x\in G$. 
\end{proof}

\begin{theorem}\label{thm: standard isometry of normed gyrogroup}
Let $G$ be a normed gyrogroup. Every transformation of the form $L_a\circ \gamma$, where $a\in G$ and $\gamma\in\GYR{G}$, is an isometry of $G$.
\end{theorem}
\begin{proof}
By Theorem 4 of \cite{TS2019MSN}, $L_a$ is an isometry of $G$ for all $a\in G$. By Corollary 2 of \cite{TS2019MSN}, $\gyr{a, b}{}$ is an isometry of $G$ for all $a, b\in G$. This implies that every element of $\GYR{G}$ is an isometry of $G$. Let $T = L_a\circ \gamma$ with $a\in G, \gamma\in\GYR{G}$. Then $T$ is an isometry of $G$, being the composite of isometries of $G$.
\end{proof}

\begin{corollary}\label{cor: gyronorm as invariant}
If $G$ is a normed gyrogroup, then the gyronorm metric $d_\eta$ is an invariant of the geometry $(G, \Gamma_m)$.
\end{corollary}

Using Corollary \ref{cor: gyronorm as invariant}, we can construct another numerical invariant of the geometry $(G, \Gamma_m)$, following the standard result in the theory of M\"{o}bius transformations. Let $G$ be a normed gyrogroup. For all points $u, v, x, y$ in $G$ such that $u\ne v$ and $x\ne y$,  define the {\it cross ratio} of $(u, v, x, y)$ by letting
\begin{equation}
[u, v; x, y] = \dfrac{d_{\eta}(u, x)d_{\eta}(v, y)}{d_{\eta}(u, v)d_{\eta}(x, y)};
\end{equation}
that is, 
\begin{equation}
[u, v; x, y] = \dfrac{\norm{\ominus u\oplus x}\norm{\ominus v\oplus y}}{\norm{\ominus u\oplus v}\norm{\ominus x\oplus y}}.
\end{equation}
By Corollary \ref{cor: gyronorm as invariant}, $[u, v; x, y] = [T(u), T(v); T(x), T(y)]$ for all  points $u, v, x, y$ in $G$ such that $u\ne v$ and $x\ne y$ and for all transformations $T$ in $\Gamma_m$. Hence, we obtain the following theorem immediately.

\begin{theorem}
If $G$ is a normed gyrogroup, then the cross ratio is an invariant of the geometry $(G, \Gamma_m)$.
\end{theorem}

In view of Theorem \ref{thm: isotropic geometry, gyrogroup with length}, every nondegenerate normed gyrogroup is {\it isotropic}, as shown in the following theorem.

\begin{theorem}
If $G$ is a nondegenerate normed gyrogroup, then the geometry $(G, \Gamma_m)$ is isotropic in the sense that for each point $p\in G$, there is a nonidentity isometry $T\in\Gamma_m$ such that $T(p) = p$.
\end{theorem}
\begin{proof}
Let $p\in G$. As in the proof of Theorem \ref{thm: isotropic geometry, gyrogroup with length}, the nonidentity transformation $T = L_p\circ \gamma\circ L_{\ominus p}$ belongs to $\Gamma_m$ and satisfies $T(p) = p$, where $\gamma$ is a nontrivial gyro-automorphism of $G$. By Theorem \ref{thm: standard isometry of normed gyrogroup}, $T$ is an isometry of $G$, which completes the proof.
\end{proof}

\begin{corollary}
If $G$ is a nondegenerate normed gyrogroup, then the geometry $(G, \Gamma_M)$ is isotropic.
\end{corollary}
\begin{proof}
The corollary follows from the fact that $\Gamma_m\subseteq \Gamma_M$.
\end{proof}

The main result of this section is presented in the following theorem, which states that the collection of open balls of equal radius is a minimally invariant set of the \mbox{geometry} $(G, \Gamma_m)$ for any normed gyrogroup $G$.

\begin{theorem}
Let $G$ be a normed gyrogroup and let $\epsilon$ be a positive number. Then $\cols{B}_{\epsilon} = \cset{B(x, \epsilon)}{x\in G}$ is a minimally invariant set of the geometry $(G, \Gamma_m)$. Here, $B(x, \epsilon)$ denotes the open ball in $G$ centered at $x$ of radius $\epsilon$; that is,
$$
B(x, \epsilon) = \cset{y\in G}{d_{\eta}(x, y) < \epsilon}.
$$
\end{theorem}
\begin{proof}
Let $a, x\in G$. By Lemma 1 of \cite{JWWATS2020}, $L_a(B(x, \epsilon)) = a\oplus B(x, \epsilon) = B(a\oplus x, \epsilon)$. Hence, $L_a(B(x, \epsilon))\in\cols{B}_{\epsilon}$. Let $\tau$ be an automorphism of $G$ that preserves the gyronorm. We claim that $\tau(B(x, \epsilon)) = B(\tau(x), \epsilon)$ for all $x\in G$. Let $y\in B(x, \epsilon)$. By Theorem 5 of \cite{TS2019MSN}, $\tau$ is an isometry of $(G, d_\eta)$. Hence, $d_\eta(\tau(x), \tau(y)) = d_\eta(x, y) < \epsilon$ and so $\tau(y)\in B(\tau(x), \epsilon)$. This shows that $\tau(B(x, \epsilon)) \subseteq B(\tau(x), \epsilon)$. Let $z\in B(\tau(x), \epsilon)$. Since $\tau$ is surjective, $z = \tau(y)$ for some $y\in G$. Note that $d_\eta(x, y) = d_\eta(\tau(x), \tau(y)) = d_\eta(\tau(x), z) < \epsilon$. Hence, $y\in B(x, \epsilon)$. This shows that $B(\tau(x), \epsilon)\subseteq \tau(B(x, \epsilon))$ and so equality holds. This in particular implies that $\gamma(B(x, \epsilon))\in\cols{B}_{\epsilon}$ for all $\gamma\in\GYR{G}$. It follows that $T(B(x, \epsilon))\in\cols{B}_{\epsilon}$ for all $T\in\Gamma_m$.

To show that $\cols{B}_{\epsilon}$ is minimally invariant, let $x, y\in G$. Then $T = L_y\circ L_{\ominus x}\in\Gamma_m$ and 
$$
T(B(x, \epsilon)) = y\oplus (\ominus x\oplus B(x, \epsilon)) = y\oplus B(e, \epsilon) = B(y, \epsilon),
$$
which completes the proof.
\end{proof}

\hspace{0.3cm}

\par\noindent\textbf{Acknowledgment.} This work was supported by the Research Center in Mathematics and Applied Mathematics, Chiang Mai University under grant No. R000025838. 

\bibliographystyle{amsplain}\addcontentsline{toc}{section}{References}
\bibliography{References}
\end{document}